\documentclass[letterpaper, 10 pt, conference]{ieeeconf}
\usepackage{amsfonts}
\usepackage{amsmath}
\usepackage{amssymb}
\usepackage{epsfig}
\usepackage{graphicx}
\usepackage[font=small,labelfont=bf,justification=justified,singlelinecheck=false]{caption}
\usepackage{subcaption}

\IEEEoverridecommandlockouts
\newtheorem{theorem}{Theorem}

\newcommand{\argmin}{\operatorname*{argmin}}
\newcommand{\ut}{\tilde{u}}

\begin{document}

\title{{\LARGE \textbf{Towards a Theoretical Analysis of PCA for
Heteroscedastic Data }}}
\author{David Hong, Laura Balzano, and Jeffrey A. Fessler
\thanks{Work by D. Hong was supported by the National Science Foundation Graduate Research Fellowship under DGE \#1256260. Work by L. Balzano was supported by the ARO Grant W911NF-14-1-0634. Work by J. Fessler was supported by the UM-SJTU data science seed fund.}}
\maketitle

\begin{abstract}
Principal Component Analysis (PCA) is a method for estimating a subspace
given noisy samples. It is useful in a variety of problems ranging from
dimensionality reduction to anomaly detection and the visualization of high
dimensional data. PCA performs well in the presence of moderate noise and even with missing data,
but is also sensitive to outliers.
PCA is also known to have a phase transition when noise is independent and identically distributed; recovery of the subspace sharply declines at a threshold noise variance.
Effective use of PCA requires a rigorous understanding of these behaviors.
This paper provides a step towards an analysis of PCA for samples with
heteroscedastic noise, that is, samples that have non-uniform noise
variances and so are no longer identically distributed. In particular, we
provide a simple asymptotic prediction of the recovery of a one-dimensional
subspace from noisy heteroscedastic samples. The prediction enables: a) easy
and efficient calculation of the asymptotic performance, and b) qualitative
reasoning to understand how PCA is impacted by heteroscedasticity (such as outliers).
\end{abstract}

\thispagestyle{empty} \pagestyle{empty}



\newcommand{\sublabel}[1]{
\renewcommand\thesubsection{\Alph{subsection}}
\addtocounter{subsection}{-1}
\refstepcounter{subsection}
\label{#1}
\renewcommand\thesubsection{\thesection.\Alph{subsection}}
}

\section{Introduction}

\noindent Given noisy measurements of points from a subspace, one may estimate the
subspace with Principal Component Analysis (PCA). Estimating a $k$%
-dimensional subspace from noisy samples $y_{1},\dots ,y_{n}\in \mathbb{R}%
^{d}$ by PCA is accomplished by solving the non-convex problem 
\begin{equation} \label{eq:cost}
\hat{U}=\argmin_{U\in \mathbb{R}^{d\times
k}:U^{T}U=I}\min_{z_{i}\in \mathbb{R}^{k}}\sum_{i=1}^{n}\Vert
y_{i}-Uz_{i}\Vert _{2}^{2},
\end{equation}%
which can be done efficiently via the singular value decomposition. PCA performs well in the presence of low to moderate noise and even performs well with missing data~\cite{wright2009rpc,chatterjee2015meb}. Furthermore, for mean zero data, representing the samples in the basis produced by PCA gives coordinates that are
uncorrelated and provide a convenient representation of the data where the
relevant factors have been decoupled.

As a result of such nice properties, PCA has been applied in myriad contexts to accomplish tasks such as dimensionality reduction, anomaly detection and the visualization of high dimensional data.
A small sample of these settings include medical imaging~\cite{ardekani1999adi}, anomaly detection on
computer networks~\cite{lakhina2004dnw} and dimensionality reduction for
classification~\cite{sharma2015and}. It has also been used to model images
taken of a scene under various illuminations~\cite{basri2003lra} as well as
measurements taken in environmental monitoring~\cite%
{papadimitriou2005spd,wagner1996sdu}, to name just a few.

To use PCA effectively in all these settings, it is important to rigorously
understand its performance under a variety of conditions. It is known, for
example, that PCA is sensitive to outliers (i.e., gross errors)~\cite{jolliffe1986pca}. Thus for problems such as
computer vision modeling~\cite{torre2001rpc} or foreground-background
separation~\cite{he2012igo} where outliers may be expected (or even of
interest), robust variants~\cite{wright2009rpc} are used instead. PCA with
independent identically distributed noise is also known to exhibit a phase
transition; recovery of the subspace sharply declines after the noise
variance exceeds a threshold~\cite{benaych2012tsv}.

This paper provides a step towards extending such analysis to
the case where noise is heteroscedastic, that is, the case where samples
have non-uniform noise variances and so are no longer identically
distributed. In particular, we provide a simple asymptotic prediction of the
recovery of a one-dimensional subspace from noisy heteroscedastic samples.
Forming the prediction involves connecting several results from random
matrix theory to obtain an initial complicated asymptotic prediction and then exploiting its structure to find a much simpler algebraic description.

The simple form enables: a) easy and efficient calculation of the asymptotic
prediction, and b) reasoning qualitatively about the expressions to
understand the asymptotic behavior of PCA with heteroscedastic noise. We
demonstrate these benefits through an example calculation and a qualitative analysis that explains a surprising phenomenon: the largest
noise variance seems to most heavily influence performance. We also perform numerical experiments to illustrate how the asymptotic prediction
applies for particular (finite) choices of ambient dimension and number of samples.

The rest of the paper is organized as follows. Section~\ref{sct:result}
describes the model we consider (a one-dimensional signal in heteroscedastic
noise) and states the main result: an asymptotic prediction for the recovery
of the one-dimensional subspace by PCA.
It also includes an example calculation of the asymptotic prediction for a particular set of model parameters, illustrating how the main result enables easy and efficient calculation of the prediction.
Section~\ref{sct:experiment}
compares the prediction with experimental results simulated according to the
model. The simulations demonstrate good agreement as the ambient dimension and number
of samples grow large; when these values are small the prediction and experiment
differ but have the same behavior.
Section~\ref{sct:proof} provides a proof of the main result.
Section~\ref{sct:analysis} uses the main
result to provide a qualitative analysis of the behavior of PCA under
heteroscedastic noise, revealing some interesting phenomena about the
negative impact of heteroscedasticity. Finally, Section~\ref{sct:discussion}
discusses the findings and describes avenues
for future work.

\vfill
\pagebreak

\section{Main result}

\label{sct:result}
\noindent We model $n$ heteroscedastic samples $y_{1},\ldots,y_{n}%
\in \mathbb{R}^{d}$ from a one dimensional subspace $\ut \in \mathbb{R}^{d}$ as%
\begin{equation} \label{eq:model}
y_{i}=\theta \ut z_{i}+\eta_{i}\varepsilon_{i}
\end{equation}
where

\begin{itemize}
\item $\theta \in \mathbb{R}_+$ is the subspace amplitude,
\item $\eta_i \in \mathbb{R}_+$ are the noise
standard deviations,

\item $\ut \in \mathbb{R}^d$ is the subspace and has entries $\ut_j \overset{\mathrm{iid}}{\sim} \mathcal{F}%
_1(0,1/d)$ with mean zero and variance $1/d$,

\item $z_i \overset{\mathrm{iid}}{\sim} \mathcal{F}_2(0,1)$ are random subspace coefficients and have mean zero and
unit variance, and

\item $\varepsilon_i\in\mathbb{R}^d$ are independent noise vectors that have entries $%
\varepsilon_{ij} \overset{\mathrm{iid}}{\sim} \mathcal{F}_3(0,1)$ with mean zero
and unit variance,
\end{itemize}

\noindent such that the distributions $\mathcal{F}_1$ and $\mathcal{F}_2$ satisfy the log-Sobolev
inequality~\cite{anderson2009ait} and the distribution $\mathcal{F}_3$ satisfies condition (1.3) from~\cite{pan2010sco}. Notably, these conditions are satisfied by Gaussian
distributions $\mathcal{F}_1=\mathcal{F}_2=\mathcal{F}_3=\mathcal{N}$.

We further suppose that $L$ noise levels $\sigma_1,\dots,\sigma_L$ occur in
proportions $p_1,\dots,p_L$. Namely, $p_1$ of the samples have noise level $%
\eta_i = \sigma_1$, $p_2$ have $\eta_i = \sigma_2$ and so on,
where the $p_\ell$ values sum to unity.

The following theorem is our main result
and describes how well the subspace $\ut$
is recovered by PCA
as the problem dimensions grow.

\medskip

\begin{theorem}
\label{thm:rs} For fixed samples-to-dimension ratio $c>\ 1$, the PCA
estimate $\hat{u}$ is such that 
\begin{equation} \label{eq:rs}
\left\vert \ut^{T}\hat{u}\right\vert ^{2}\ \ \underset{{\tiny 
\begin{array}{c}
n,d\rightarrow \infty \\ 
n/d=c%
\end{array}%
}}{\overset{a.s.}{\longrightarrow }}\max \left( 0,\frac{A\left( \beta
\right) }{\beta B^{\prime }\left( \beta \right) }\right)
\end{equation}%
where%
\begin{align*}
A\left( x\right) & =1-c\sum_{\ell =1}^{L}\frac{p_{\ell }\sigma _{\ell }^{4}}{%
\left( x-\sigma _{\ell }^{2}\right) ^{2}} \\
B\left( x\right) & =1-c\theta ^{2}\sum_{\ell =1}^{L}\frac{p_{\ell }}{%
x-\sigma _{\ell }^{2}}
\end{align*}%
and $\beta $ is the largest real root of $B$.
\end{theorem}

\medskip

\noindent Section~\ref{sct:proof} presents the proof of this theorem. We illustrate Theorem~\ref{thm:rs} with the following example calculation.
\\


\noindent \textbf{Example calculation:} Here we calculate the asymptotic prediction in~\eqref{eq:rs} for the case: 
\begin{align*}
c&=5 & p&=\left(0.2,0.8\right) \\
\theta&=2 & \sigma&=\left(1,2\right)
\end{align*}
Namely, we determine the limit of $\left\vert \ut^{T}\hat{u}\right\vert ^{2}$
when $d,n\rightarrow\infty$ for the case where there are $5$ times as many
samples as the ambient dimension, the signal amplitude is $2$, and $20\%$ of
the samples have low noise with variance $1$ (signal-to-noise ratio $%
\theta^{2}/\sigma_1^{2} = 4$) and $80\%$ of the samples have high noise with
variance $4$ (signal-to-noise ratio $\theta^{2}/\sigma_2^{2} = 1$).
The steps are as follows.
\begin{enumerate}
\item Substitute the values of $c,\theta,p$ and $\sigma$ into the formulas
for $A$ and $B$, obtaining 
\begin{align*}
A\left( x\right) & =1-5\cdot\left( \frac{0.2\ \cdot1^{4}}{\left(
x-1^{2}\right) ^{2}}+\frac{0.8\ \cdot2^{4}}{\left( x-2^{2}\right) ^{2}}%
\right) \\
& =1-\frac{1}{\left( x-1\right) ^{2}}-\frac{64}{\left( x-4\right) ^{2}} \\
B\left( x\right) & =1-5\cdot2^{2}\cdot\left( \frac{0.2}{x-1^{2}}+\frac{0.8}{%
x-2^{2}}\right) \\
& =1-\frac{4}{x-1}-\frac{16}{x-4}.
\end{align*}

\item Find the largest root of $B$, obtaining 
\begin{equation*}
\beta =23.466.
\end{equation*}

\item Evaluate $\frac{A\left( \beta \right) }{\beta B^{\prime }\left( \beta
\right) },$ obtaining%
\begin{equation*}
\frac{A\left( \beta \right) }{\beta B^{\prime }\left( \beta \right) }=0.705.
\end{equation*}

\item Take the maximum with zero, and conclude that%
\begin{equation*}
\left\vert \ut^{T}\hat{u}\right\vert ^{2}\ \ \underset{{\tiny 
\begin{array}{c}
n,d\rightarrow \infty \\ 
n/d=c%
\end{array}%
}}{\overset{a.s.}{\longrightarrow }}\ \ 0.705.
\end{equation*}
\end{enumerate}

\noindent Note that the second step can be easily done by clearing the
denominator of $B$ and finding the real roots of the resulting degree $L$
polynomial (using off-the-shelf tools).
Hence the asymptotic prediction can be efficiently computed.

\section{Experimental verification}

\label{sct:experiment}

\noindent To illustrate the main result (Theorem~\ref{thm:rs})
we performed a numerical experiment for the two noise level case ($L=2$):
\begin{align*}
c &= 10 & \theta&=1 & \sigma&=(1.8,0.2)
\end{align*}
where $p_2$ is swept from $0$ to $1$ with $p_1 = 1-p_2$. This allows us to
investigate the accuracy of the asymptotic prediction for a variety of
settings: at the extremes ($p_2=0,1$) the setup matches the
homoscedastic setting and in the middle ($p_2=1/2$) the samples are
split evenly between the two noise levels.

\begin{figure}[t]
\hrule
\centering
\begin{subfigure}[t]{0.47\linewidth}
\includegraphics[width=\linewidth]{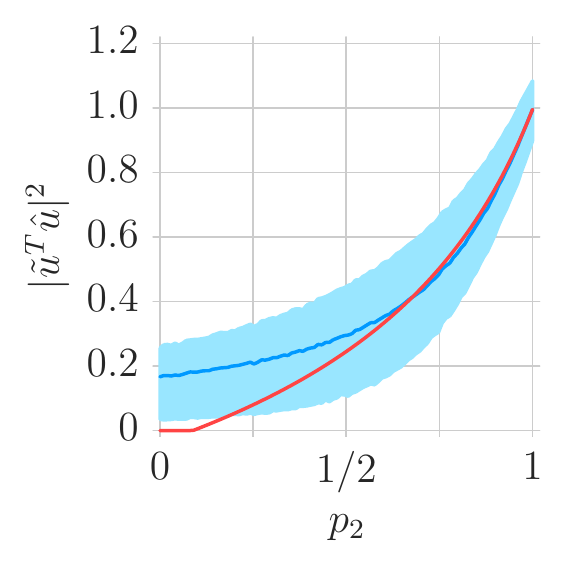}
\caption{Results for $d=10^2$ and $n=10^3$ ($%
10000$ trials).}
\label{fig:exp1}
\end{subfigure}
\hfill
\begin{subfigure}[t]{0.47\linewidth}
\includegraphics[width=\linewidth]{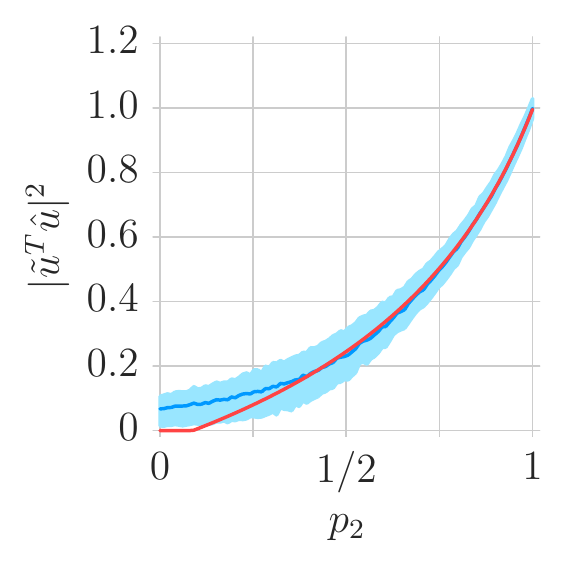}
\caption{Results for $d=10^3$ and $n=10^4$ ($%
1000 $ trials).}
\label{fig:exp2}
\end{subfigure}
\caption{Simulation results for $c = 10$, $\theta=1$, $\sigma=(1.8,0.2)$
where $p_2$ is swept from $0$ to $1$ with $p_1 = 1-p_2$. Simulation mean (blue curve) and interquartile
interval (light blue ribbon) shown with asymptotic prediction (red curve).}
\label{fig:exp}
\hrule
\end{figure}

We first suppose $d=10^2$ and $n=10^3$. We performed $10000$ trials with
data generated by the Gaussian distribution. Namely, 
$\mathcal{F}_1=\mathcal{F}_2=\mathcal{F}_3=\mathcal{N}$.
Figure~\ref{fig:exp1} shows the simulation results with the mean (blue
curve) and interquartile interval (light blue ribbon) shown with the
asymptotic prediction (red curve).

There is generally good agreement between the mean and the
asymptotic prediction for $p_2 > 0.6$ (they deviate from each other by no
more than $0.025$). However, the smaller the value of~$p_2$, the greater the deviation
of the prediction from the mean, with the asymptotic prediction
underestimating the (non-asymptotic) simulation by at most~$0.182$. 

Figure~\ref{fig:exp} illustrates a general phenomenon we observed: small asymptotic predictions typically
underestimate the simulation results, with smaller predictions
underestimating more. Intuitively, a small prediction of this inner product
 corresponds to a subspace estimate that
is becoming increasingly like an isotropically random vector and
so it has vanishing square inner product with the true subspace as the
dimension grows. 
Asymptotically, below the
phase transition, theory predicts the subspace estimate to be practically isotropically random.
However, in finite dimension there is a better chance of alignment, resulting in 
a positive square inner product.

We illustrate this phenomenon with a second experiment with a higher
dimension of $d=10^3$ and a higher number of samples $n=10^4$ (chosen to
have the same sample-to-dimension ratio $c=10$). The data were again generated
using the Gaussian distribution and we performed $1000$ trials. Figure~\ref%
{fig:exp2} shows the simulation results for this case. Again, the mean and
interquartile interval are in blue and light blue, and the asymptotic prediction is in
red. This experiment demonstrates better agreement between the mean behavior and the asymptotic
prediction. In particular, for $p_2 > 0.4$ they deviate from each other by
no more than $0.017$. For $p_2 < 0.4$ their largest deviation is $0.084$;
notably, this is less than half of the largest deviation for the smaller
experiment. Furthermore, the interquartile interval is narrower, indicating
that the square inner product is concentrating.

We stress here that, while we have not proved this relationship, 
experimental evidence suggests that the asymptotic 
prediction in Theorem~\ref{thm:rs} underestimates the mean square inner product obtained in
finite size experiments, and is therefore a conservative or pessimistic estimate
of the performance of PCA. When determining how much to trust the subspace estimated by PCA, this conservatism might be preferable to overestimating the
reliability of PCA.

Finally, note that the square inner product from the simulation is sometimes
above one. This is because the true $\ut$ is a random vector and may not have
unit norm. However, as the dimension grows, the norm of $\ut$ will concentrate
around one (see~\cite{vincent2010com} for a treatment of the
concentration of the norm).

\section{Proof of main result}
\label{sct:proof}
\noindent This section proves the main result (Theorem~\ref{thm:rs}). The proof
has eight main parts. In~\ref{sct:initial}, we first apply previous results
in the literature to obtain an initial expression for the asymptotic
prediction. This expression is difficult to evaluate and analyze because it
involves an integral transform of the (nontrivial) limiting singular value distribution
for a random (noise) matrix as well as the corresponding limiting largest
singular value. In the remaining seven parts (\ref{sct:varchange}-\ref{sct:r_max}), we find a simple equivalent expression by exploiting the structure of the prediction. 

\subsection{Obtain an initial expression}

\label{sct:initial}%

\noindent Rewriting the model in~\eqref{eq:model} in matrix form yields%
\begin{equation*}
\mathbf{Y}:=\theta \ut z^{T}+\mathbf{E} \mathbf{H} \in \mathbb{R}^{d\times n}
\end{equation*}%
where $\mathbf{E} \in \mathbb{R}^{d\times n}$ is a matrix with columns $%
\varepsilon _{1},\ldots ,\varepsilon _{n}\in \mathbb{R}^{d}$ and $\mathbf{H} \in 
\mathbb{R}^{n\times n}$ is a diagonal matrix with diagonal entries $\eta
_{1},\ldots ,\eta _{n}\in \mathbb{R}$.

Recall that the subspace basis $\hat{u}$ estimated by PCA for $\mathbf{Y}$ is the
first left singular vector of $\mathbf{Y}$. PCA is invariant to scaling, so $\hat{u}$ is also the first left singular vector of%
\begin{equation*}
\mathbf{\tilde{Y}}:=\frac{1}{\sqrt{n}}\mathbf{Y}.
\end{equation*}
The matrix $\mathbf{\tilde{Y}}$ matches the low rank (here rank one) perturbation
of a random matrix model considered in\ \cite{benaych2012tsv} because%
\begin{equation*}
\mathbf{\tilde{Y}}=\mathbf{P}+\mathbf{X}
\end{equation*}%
where%
\begin{equation*}
\mathbf{P}:=\theta \ut\left( \frac{1}{\sqrt{n}}z\right) ^{T}\mathbf{X}:=\left( \frac{1}{\sqrt{n}}%
\mathbf{E} \right) \mathbf{H} ,
\end{equation*}%
and $\mathbf{P}$ is generated according to the \textquotedblleft i.i.d
model\textquotedblright\ and satisfies Assumption 2.4 of~\cite%
{benaych2012tsv},\ and $\mathbf{X}$ satisfies Assumptions 2.1-2.3 of~\cite%
{benaych2012tsv} ($\mathbf{X}$ here matches the random matrix in~\cite{pan2010sco},
which~\cite{benaych2012tsv} refers to as an example of a
random matrix that satisfies the assumptions).

Thus under the condition $\varphi ^{\prime }\left( b^{+}\right) =-\infty $ (we will show in
subsection~\ref{sct:props} that it is indeed satisfied), Theorems 2.10 and 2.11 from~\cite{benaych2012tsv} yield%
\begin{equation} \label{eq:init}
\left\vert \ut^{T}\hat{u}\right\vert ^{2}\ \ \underset{{\tiny 
\begin{array}{c}
n,d\rightarrow \infty \\ 
n/d=c%
\end{array}%
}}{\overset{a.s.}{\longrightarrow }}\left\{ 
\begin{array}{ll}
\frac{-2\varphi \left( \rho \right) }{\theta ^{2}D^{\prime }\left( \rho
\right) } & \theta ^{2}>\bar{\theta}^{2} \\ 
0 & \text{otherwise}%
\end{array}%
\right.
\end{equation}%
where

\begin{itemize}
\item $\rho :=D^{-1}\left( 1/\theta ^{2}\right) $

\item $\bar{\theta}^{2}:=1/D\left( b^{+}\right) $

\item $D\left( z\right) :=\varphi \left( z\right) \left( c^{-1}\varphi
\left( z\right) +\frac{1-c^{-1}}{z}\right) \qquad z>b$

\item $\varphi \left( z\right) :=\int_{a}^{b}\frac{z}{z^{2}-t^{2}}d\mu
_{\mathbf{X}}\left( t\right) \qquad z>b$

\item $a$ and $b$ are, respectively, the infimum and the supremum of the
support of $\mu _{\mathbf{X}}$ (so $b >a\geq0$), and

\item $\mu _{\mathbf{X}}$ is the limiting singular value distribution of $\mathbf{X}$
(compactly supported by Assumption 2.1 of~\cite{benaych2012tsv}).
\end{itemize}

\noindent We use the notation $f\left( b^{+}\right) :=\lim_{z\rightarrow b^{+}}f\left( z\right) $ as a
convenient shorthand for the limit of a function $f\left( z\right)$.

Evaluating this asymptotic prediction would then consist of evaluating the
above intermediates from bottom to top. These steps are challenging because
because they
involve an integral transform of the limiting singular value distribution
for the random (noise) matrix as well as the corresponding limiting largest
singular value.
 The following sections eventually lead to the
simpler expression in~\eqref{eq:rs} that is easier to both evaluate and analyze.

\subsection{Carry out a change of variables}

\label{sct:varchange}%

\noindent We begin by introducing the function%
\begin{equation} \label{eq:psi}
\psi \left( z\right) :=\frac{cz}{\varphi \left( z\right) }=\left[ \frac{1}{c}%
\int_{a}^{b}\frac{1}{z^{2}-t^{2}}d\mu _{\mathbf{X}}\left( t\right) \right]
^{-1},\ z>b
\end{equation}%
because it turns out to have several nice properties that simplify all of the following analysis.

Rewriting~\eqref{eq:init} using $\psi
\left( z\right) $ instead of $\varphi \left( z\right) $ yields%
\begin{equation*}
\left\vert \ut^{T}\hat{u}\right\vert ^{2}\ \ \underset{{\tiny 
\begin{array}{c}
n,d\rightarrow \infty \\ 
n/d=c%
\end{array}%
}}{\overset{a.s.}{\longrightarrow }}\left\{ 
\begin{array}{ll}
\frac{-2c}{\theta ^{2}\psi \left( \rho \right) D^{\prime }\left( \rho
\right) /\rho } & \theta ^{2}>\bar{\theta}^{2} \\ 
0 & \text{otherwise}%
\end{array}%
\right.
\end{equation*}%
where now
\begin{equation*}
D\left( z\right) =\frac{cz^{2}}{\left( \psi \left( z\right) \right)
^{2}}+\frac{c-1}{\psi \left( z\right) },\quad z>b.
\end{equation*}

\subsection{Find useful properties of $\protect\psi \left( z\right) $}

\label{sct:props}%

\noindent Establishing some properties of $\psi \left( z\right) 
$ aids simplification significantly. Furthermore, these properties help us
show that $\varphi ^{\prime }\left( b^{+}\right) $ is indeed $-\infty $, as
stated above in subsection~\ref{sct:initial}.

\medskip

\noindent \textbf{Property 1.} We show that $\psi \left( z\right) $
satisfies a certain rational equation for all $z>b$. For this, observe that
the square singular values of the noise matrix $\mathbf{X}$ are exactly the eigenvalues of $c\mathbf{X}\mathbf{X}^{T}$
divided by $c$ (since $\mathbf{X}$ has more columns than rows, namely, $c=n/d>1$).
Thus we first consider the limiting eigenvalue distribution $\mu _{c\mathbf{X}\mathbf{X}^{T}}$
of $c\mathbf{X}\mathbf{X}^{T}$, and then relate its Stieltjes transform $m\left( \zeta \right) 
$ to $\psi \left( z\right) $.

Theorem~1 in~\cite{pan2010sco} establishes that the random matrix%
\begin{equation*}
c\mathbf{X}\mathbf{X}^{T}=\left( \frac{1}{\sqrt{d}}\mathbf{E} \right) \mathbf{H} ^{2}\left( \frac{1%
}{\sqrt{d}}\mathbf{E} \right) ^{T}
\end{equation*}%
has a limiting eigenvalue distribution $\mu _{c\mathbf{X}\mathbf{X}^{T}}$ whose Stieltjes
transform%
\begin{equation} \label{eq:stieltjes}
m\left( \zeta \right) :=\int \frac{1}{t-\zeta }d\mu _{c\mathbf{X}\mathbf{X}^{T}}\left(
t\right) ,\quad \zeta \in \mathbb{C}^{+}
\end{equation}%
satisfies the condition%
\begin{equation} \label{eq:pan}
\forall \zeta \in \mathbb{C}^{+}\quad m\left( \zeta \right) =-\left( \zeta
-c\sum_{\ell =1}^{L}\frac{p_{\ell }\sigma _{\ell }^{2}}{1+\sigma _{\ell
}^{2}m\left( \zeta \right) }\right) ^{-1}
\end{equation}%
where $\mathbb{C}^{+}$ is the set of all complex numbers with positive
imaginary part.

Since the square singular values of $\mathbf{X}$ are exactly the eigenvalues of $%
c\mathbf{X}\mathbf{X}^{T}$ divided by $c$, we have for all $z > b$%
\begin{align}
\psi \left( z\right) &=\left[ \frac{1}{c}\int_{a}^{b}\frac{1}{z^{2}-t^{2}}%
d\mu _{\mathbf{X}}\left( t\right) \right] ^{-1} \nonumber\\
&=\left[ \frac{1}{c}\int_{a^{2}c}^{b^{2}c}\frac{1}{z^{2}-t/c}d\mu
_{c\mathbf{X}\mathbf{X}^{T}}\left( t\right) \right] ^{-1} \nonumber\\
&=-\left[ \int_{a^{2}c}^{b^{2}c}\frac{1}{t-z^{2}c}d\mu _{c\mathbf{X}\mathbf{X}^{T}}\left(
t\right) \right] ^{-1}. \label{eq:X_cXXt}
\end{align}
For all $z$ and $\xi >0$, $z^{2}c+i\xi \in \mathbb{C}^{+}$
and so combining~\eqref{eq:stieltjes}-\eqref{eq:X_cXXt} yields that for all $z>b$%
\begin{align*}
\psi \left( z\right) &=-\left[ \lim_{\xi \rightarrow 0^{+}}m\left(
z^{2}c+i\xi \right) \right] ^{-1} \\
&=\lim_{\xi \rightarrow 0^{+}}z^{2}c+i\xi -c\sum_{\ell =1}^{L}\frac{p_{\ell
}\sigma _{\ell }^{2}}{1+\sigma _{\ell }^{2}m\left( z^{2}c+i\xi \right) } \\
&=z^{2}c -c\sum_{\ell =1}^{L}\frac{p_{\ell }\sigma _{\ell }^{2}}{1+\sigma
_{\ell }^{2}\lim_{\xi \rightarrow 0^{+}}m\left( z^{2}c+i\xi \right) } \\
&=z^{2}c-c\sum_{\ell =1}^{L}\frac{p_{\ell }\sigma _{\ell }^{2}}{1-\sigma
_{\ell }^{2}/\psi \left( z\right) }.
\end{align*}
Rearranging yields%
\begin{equation} \label{eq:Q1}
\forall z>b,\  0=\frac{cz^{2}}{\left( \psi \left( z\right) \right) ^{2}}-%
\frac{1}{\psi \left( z\right) }-\frac{c}{\psi \left( z\right) }\sum_{\ell
=1}^{L}\frac{p_{\ell }\sigma _{\ell }^{2}}{\psi \left( z\right) -\sigma
_{\ell }^{2}}.
\end{equation}
where the last term is
\begin{align*}
-\frac{c}{\psi \left( z\right) }\sum_{\ell =1}^{L}\frac{p_{\ell }\sigma
_{\ell }^{2}}{\psi \left( z\right) -\sigma _{\ell }^{2}} &=\frac{c}{\psi
\left( z\right) }\sum_{\ell =1}^{L}p_{\ell }\left[ 1-\frac{\psi \left(
z\right) }{\psi \left( z\right) -\sigma _{\ell }^{2}}\right] \\
&=\frac{c}{\psi \left( z\right) }-c\sum_{\ell =1}^{L}\frac{p_{\ell }}{\psi
\left( z\right) -\sigma _{\ell }^{2}}
\end{align*}%
because $p_{1}+\cdots +p_{L}=1$. Substituting back into~\eqref{eq:Q1} finally yields $%
0=Q\left( \psi \left( z\right) ,z\right) $ for all $z>b$, where%
\begin{equation} \label{eq:Q}
Q\left( s,z\right) :=\frac{cz^{2}}{s^{2}}+\frac{c-1}{s}-c\sum_{\ell =1}^{L}%
\frac{p_{\ell }}{s-\sigma _{\ell }^{2}}
.
\end{equation}%
Thus $\psi $ is an algebraic function with associated rational function $Q$
(a polynomial can be formed by clearing the denominator).

\medskip

\noindent \textbf{Property 2. }We show that $\psi \left( b^{+}\right) $ is
finite and $\psi ^{\prime }\left( b^{+}\right) =\infty $. For this, note
first that $\psi \left( b^{+}\right) $ is a multiple root of $Q\left(
\cdot ,b\right)$ and hence is finite. This follows
from the observation in \cite{nadakuditi2006tpm} that non-pole boundary
points of compactly supported distributions like $\mu _{c\mathbf{X}\mathbf{X}^{T}}$ occur
where the polynomial defining the Stieltjes transform has multiple roots.

Differentiating $0=Q\left( \psi \left( z\right) ,z\right) $ with respect
to $z$ and rearranging yields%
\begin{equation*}
\psi ^{\prime }\left( z\right) =-\frac{\frac{\partial Q}{\partial z}\left(
\psi \left( z\right) ,z\right) }{\frac{\partial Q}{\partial s}\left( \psi
\left( z\right) ,z\right) }.
\end{equation*}%
Since $\psi \left( b^{+}\right) $ is a multiple root of $Q\left( \cdot
,b\right) $,%
\begin{equation*}
\frac{\partial Q}{\partial s}\left( \psi \left( b^{+}\right) ,b\right) =0
\end{equation*}%
while on the other hand%
\begin{equation*}
\frac{\partial Q}{\partial z}\left( \psi \left( b^{+}\right) ,b\right) =%
\frac{2cb}{\left( \psi \left( b^{+}\right) \right) ^{2}}>0.
\end{equation*}%
Thus $\psi ^{\prime }\left( b^{+}\right) =\infty $, where
the sign is necessarily positive because $\psi \left( z\right) $ is an
increasing function.

\medskip
\newpage
\noindent Summarizing, we have shown that

\begin{enumerate}
\item $\psi $ satisfies the equation $0=Q\left( \psi \left( z\right)
,z\right) $ for all $z>b$

\item $\psi \left( b^{+}\right) $ is finite and $\psi ^{\prime }\left(
b^{+}\right) =\infty $
\end{enumerate}

\noindent As an immediate consequence of these properties, we also have that
indeed

\begin{equation*}
\varphi ^{\prime }\left( b^{+}\right) =\frac{c}{\psi \left( b^{+}\right) }%
\left[ 1-z\frac{\psi ^{\prime }\left( b^{+}\right) }{\psi \left(
b^{+}\right) }\right] =-\infty .
\end{equation*}

\subsection{Express $D\left( z\right) $ in terms of only $\protect\psi %
\left( z\right) $}

\label{sct:D_psi}%

\noindent This subsection uses the properties of $\psi \left( z\right) 
$ to find a simple expression for $D\left( z\right) $ in terms of $\psi
\left( z\right) $. Observe that%
\begin{equation*}
D\left( z\right) =Q\left( \psi \left( z\right) ,z\right) +c\sum_{\ell =1}^{L}%
\frac{p_{\ell }}{\psi \left( z\right) -\sigma _{\ell}^{2}} .
\end{equation*}%
Recalling that $0=Q\left( \psi \left( z\right) ,z\right) $ for $z>b$, we have%
\begin{equation} \label{eq:D}
D\left( z\right) =c\sum_{\ell =1}^{L}\frac{p_{\ell }}{\psi \left( z\right)
-\sigma _{\ell }^{2}} .
\end{equation}

\subsection{Express $D^{\prime }\left( z\right) /z$ in terms of only $%
\protect\psi \left( z\right) $}

\label{sct:Dp_psi}%

\noindent This subsection uses the properties of $\psi \left( z\right) 
$ to find a simple expression for $D^{\prime }\left( z\right) /z$ in terms
of $\psi \left( z\right) $. Differentiating~\eqref{eq:D}
with respect to $z$ yields%
\begin{equation*}
D^{\prime }\left( z\right) =-c\psi ^{\prime }\left( z\right) \sum_{\ell
=1}^{L}\frac{p_{\ell }}{\left( \psi \left( z\right) -\sigma _{\ell
}^{2}\right) ^{2}}
\end{equation*}%
and so we need to express $\psi ^{\prime }\left( z\right) $ in terms of $%
\psi \left( z\right) $.

To do this, differentiate both sides of $0=Q\left( \psi \left( z\right)
,z\right) $ with respect to $z$ and solve for $\psi ^{\prime }\left(
z\right) $, obtaining%
\begin{equation*}
\psi ^{\prime }\left( z\right) =\frac{2cz}{\gamma \left( z\right) }
\end{equation*}%
where the denominator is%
\begin{equation*}
\gamma \left( z\right) :=c-1+\frac{2cz^{2}}{\psi \left( z\right) }%
-c\sum_{\ell =1}^{L}\frac{p_{\ell }\left( \psi \left( z\right) \right) ^{2}}{%
\left( \psi \left( z\right) -\sigma _{\ell }^{2}\right) ^{2}} .
\end{equation*}
Note that%
\begin{equation*}
\frac{2cz^{2}}{\psi \left( z\right) }=-2\left( c-1\right) +c\sum_{\ell
=1}^{L}\frac{2p_{\ell }\psi \left( z\right) }{\psi \left( z\right) -\sigma
_{\ell }^{2}}
\end{equation*}%
because $0=Q\left( \psi \left( z\right) ,z\right) $ for $z>b$. Substituting into $\gamma(z)$ and forming a common denominator yields%
\begin{align*}
\gamma \left( z\right) &=1-c +c\sum_{\ell =1}^{L}\frac{2p_{\ell }\psi \left(
z\right) }{\psi \left( z\right) -\sigma _{\ell }^{2}}-c\sum_{\ell =1}^{L}%
\frac{p_{\ell }\left( \psi \left( z\right) \right) ^{2}}{\left( \psi \left(
z\right) -\sigma _{\ell }^{2}\right) ^{2}} \\
&=1-c +c\sum_{\ell =1}^{L}p_{\ell }\frac{\left( \psi \left(z\right) \right)
^{2}-2\psi \left( z\right) \sigma _{\ell }^{2}}{\left( \psi \left( z\right)
-\sigma _{\ell }^{2}\right) ^{2}}
\end{align*}
Dividing the summand with respect to $\psi(z)$ and recalling that $%
p_{1}+\cdots +p_{L}=1$ yields
\begin{align*}
\gamma(z) &=1-c +c\sum_{\ell =1}^{L}\left( p_{\ell }-\frac{p_{\ell}\sigma
_{\ell }^{4}}{\left( \psi \left( z\right) -\sigma _{\ell}^{2}\right) ^{2}}%
\right) \\
&=1-c\sum_{\ell =1}^{L}\frac{p_{\ell }\sigma _{\ell }^{4}}{\left( \psi\left(
z\right) -\sigma _{\ell }^{2}\right) ^{2}} =A\left( \psi \left(z\right)
\right)
\end{align*}%
where%
\begin{equation*}
A\left( x\right) :=1-c\sum_{\ell =1}^{L}\frac{p_{\ell }\sigma _{\ell }^{4}}{%
\left( x-\sigma _{\ell }^{2}\right) ^{2}}.
\end{equation*}
Thus%
\begin{equation} \label{eq:psip}
\psi ^{\prime }\left( z\right) =\frac{2cz}{A\left( \psi \left( z\right)
\right) }
\end{equation}
and so%
\begin{equation} \label{eq:Dp}
\frac{D^{\prime }\left( z\right) }{z}=-\frac{2c^{2}}{A\left( \psi \left(
z\right) \right) }\sum_{\ell =1}^{L}\frac{p_{\ell }}{\left( \psi \left(
z\right) -\sigma _{\ell }^{2}\right) ^{2}}.
\end{equation}

\subsection{Express the prediction in terms of only $\protect\psi \left(
b^{+}\right) $ and $\protect\psi \left( \protect\rho \right) $}

\label{sct:r_psi}%

\noindent This subsection uses~\eqref{eq:D} and~\eqref{eq:Dp} to express the asymptotic prediction in terms
of $\psi \left( b^{+}\right) $ and $\psi \left( \rho \right) $. Using~\eqref{eq:D} yields
\begin{equation*}
\frac{1}{\bar{\theta}^{2}}=D\left( b^{+}\right) =c\sum_{\ell =1}^{L}\frac{%
p_{\ell }}{\psi \left( b^{+}\right) -\sigma _{\ell }^{2}}.
\end{equation*}%
Thus the condition $\theta ^{2}>\bar{\theta}^{2}$ is equivalent to%
\begin{equation*}
0>1-\frac{\theta ^{2}}{\bar{\theta}^{2}}=1-c\theta ^{2}\sum_{\ell =1}^{L}%
\frac{p_{\ell }}{\psi \left( b^{+}\right) -\sigma _{\ell }^{2}}=B\left( \psi
\left( b^{+}\right) \right)
\end{equation*}%
where%
\begin{equation*}
B\left( x\right) :=1-c\theta ^{2}\sum_{\ell =1}^{L}\frac{p_{\ell }}{x-\sigma
_{\ell }^{2}}.
\end{equation*}
Using~\eqref{eq:Dp} yields
\begin{align*}
r&:=\frac{-2c}{\theta ^{2}\psi \left( \rho \right) D^{\prime }\left( \rho
\right) /\rho } =\frac{A\left( \psi \left( \rho \right) \right) }{\psi
\left( \rho \right) c\theta ^{2}\sum_{\ell =1}^{L}\frac{p_{\ell }}{\left(
\psi \left( \rho \right) -\sigma _{\ell }^{2}\right) ^{2}}} \\
&=\frac{A\left( \psi \left( \rho \right) \right) }{\psi \left( \rho \right)
B^{\prime }\left( \psi \left( \rho \right) \right) }
\end{align*}%
where we note that%
\begin{equation*}
B^{\prime }\left( x\right) =c\theta ^{2}\sum_{\ell =1}^{L}\frac{p_{\ell }}{%
\left( x-\sigma _{\ell }^{2}\right) ^{2}}.
\end{equation*}
Summarizing, the asymptotic prediction is now expressed as%
\begin{equation*}
\left\vert \ut^{T}\hat{u}\right\vert ^{2}\ \ \underset{{\tiny 
\begin{array}{c}
n,d\rightarrow \infty \\ 
n/d=c%
\end{array}%
}}{\overset{a.s.}{\longrightarrow }}\left\{ 
\begin{array}{ll}
\frac{A\left( \psi \left( \rho \right) \right) }{\psi \left( \rho \right)
B^{\prime }\left( \psi \left( \rho \right) \right) } & B\left( \psi \left(
b^{+}\right) \right) <0 \\ 
0 & \text{otherwise.}%
\end{array}%
\right. 
\end{equation*}

\subsection{Express the prediction algebraically}

\label{sct:r_alg}%

\noindent This subsection finds an algebraic description of the asymptotic
prediction. We first use the properties of $\psi \left( z\right) $
to show that $\psi \left( b^{+}\right) $ and $\psi \left( \rho \right) $
are, respectively, the largest real roots of $A$ and $B$.

Using~\eqref{eq:psip} yields
\begin{equation*}
A\left( \psi \left( b^{+}\right) \right) =\frac{2cz}{\psi ^{\prime }\left(
b^{+}\right) }=0
\end{equation*}%
because $\psi ^{\prime }\left( b^{+}\right) =\infty $. Thus $\psi \left(
b^{+}\right) $ is a real root of $A$.

For $\theta >\bar{\theta}$, we have $\rho :=D^{-1}\left( 1/\theta ^{2}\right)$
and so%
\begin{equation*}
0=1-\theta ^{2}D\left( \rho \right) =1-c\theta ^{2}\sum_{\ell =1}^{L}\frac{%
p_{\ell }}{\psi \left( \rho \right) -\sigma _{\ell }^{2}}=B\left( \psi
\left( \rho \right) \right) .
\end{equation*}%
Thus $\psi \left( \rho \right) $ is a real root of $B$.

The functions $A$ and $B$ both have several real roots. To show
that $\psi \left( b^{+}\right) $ and $\psi \left( \rho \right) $ are the
largest ones, consider%
\begin{equation*}
\psi \left( \rho \right) =\psi \left( D^{-1}\left( 1/\theta ^{2}\right)
\right)
\end{equation*}
as a function of $\theta $ as $\theta $ increases from $\bar{\theta}$ to
infinity. Note that%
\begin{equation*}
\psi \left( z\right) =\left[ \frac{1}{c}\int_{a}^{b}\frac{1}{z^{2}-t^{2}}%
d\mu _{\mathbf{X}}\left( t\right) \right]^{-1},\quad z>b
\end{equation*}%
continuously and monotonically increases from $\psi \left( b^{+}\right) $ to
infinity as $z$ increases from $b$ to infinity. Thus%
\begin{equation*}
D\left( z\right) =c\sum_{\ell =1}^{L}\frac{p_{\ell }}{\psi \left( z\right)
-\sigma _{\ell }^{2}},\quad z>b
\end{equation*}%
continuously and monotonically decreases from $1/\bar{\theta}^{2}$ to zero
as $z$ increases from $b$ to infinity, and so $D^{-1}\left( 1/\theta
^{2}\right) $ continuously and monotonically increases from $b$ to infinity
as $\theta $ increases from $\bar{\theta}$ to infinity.

\begin{figure}[t!]
\hrule
\begin{center}
\vspace{1.8mm}
\includegraphics[width=0.98\linewidth]{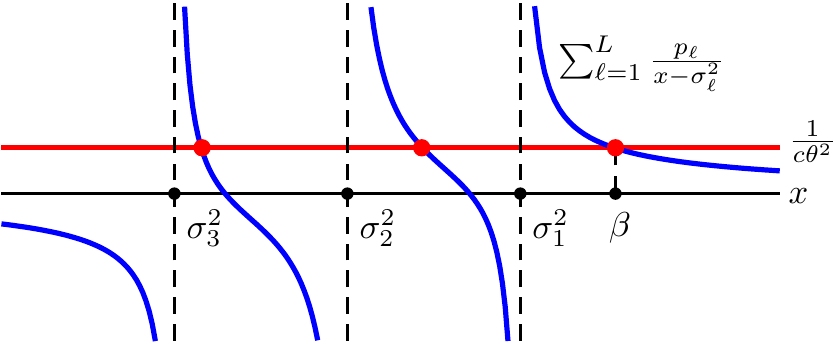}
\end{center}
\caption{Illustration of the real roots of $B(x)$ for $L=3$ levels. They
occur where the sum (blue curve) intersects $1/(c\protect\theta^2)$ (red
line). The largest is $\protect\beta$.}
\label{fig:B}
\hrule
\end{figure}

As a result, $\psi \left( \rho \right) $ continuously and monotonically
increases from $\psi \left( b^{+}\right) $ towards infinity as $\theta $
increases from $\bar{\theta}$ towards infinity. This is possible only if both $%
\psi \left( b^{+}\right) $ and $\psi \left( \rho \right) $ are larger than
all the noise levels $\sigma _{\ell }^{2}$. To see this, recall that $\psi
\left( \rho \right) $ is a real root of $B$ and note that the real roots of $B$ satisfy the following equation (illustrated in Figure~\ref{fig:B}):
\begin{equation*}
\frac{1}{c\theta ^{2}}=\sum_{\ell =1}^{L}\frac{p_{\ell }}{x-\sigma _{\ell
}^{2}}.
\end{equation*}
If either $\psi \left(
b^{+}\right) $ or $\psi \left( \rho \right) $ were less than any of the
noise levels~$\sigma _{\ell }^{2}$,
then $\psi(\rho)$ would change discontinuously as $\theta$ varies.
Thus $\psi \left( b^{+}\right) $ and $\psi \left( \rho \right) $ are indeed both larger than all the noise levels.

To the right of all the noise levels (i.e., for $x$ larger than all the
noise levels $\sigma _{\ell }^{2}$), both $A$ and $B$ continuously and
monotonically increase from negative infinity to one, and so each has
exactly one real root larger than all the noise levels (namely, the largest
real root). Thus $\psi \left( b^{+}\right) $ is the largest
real root of $A$ and when $\theta >\bar{\theta}$, $\psi \left( \rho
\right) $ is the largest real root of~$B$.

Using this yields the algebraic form of the
prediction:%
\begin{equation*}
\left\vert \ut^{T}\hat{u}\right\vert ^{2}\ \ \underset{{\tiny 
\begin{array}{c}
n,d\rightarrow \infty \\ 
n/d=c%
\end{array}%
}}{\overset{a.s.}{\longrightarrow }}\left\{ 
\begin{array}{ll}
\frac{A\left( \beta \right) }{\beta B^{\prime }\left( \beta \right) } & 
B\left( \alpha \right) <0 \\ 
0 & \text{otherwise}%
\end{array}%
\right.
\end{equation*}%
where $\alpha $ and $\beta $ are, respectively, the largest real roots of $A$
and $B$.

\subsection{Further simplify the asymptotic prediction}

\label{sct:r_max}%

\noindent We further simplify the asymptotic prediction by showing that $%
B\left( \alpha \right) <0$ is equivalent to $A\left( \beta \right) /\left(
\beta B^{\prime }\left( \beta \right) \right) >0$.

To do this, observe that both $\alpha $ and $\beta $ are larger than all the
noise levels $\sigma _{\ell }^{2}$, and note that $A\left( x\right) $ and $%
B\left( x\right) $ are both monotonically increasing in this regime. Thus it
follows that%
\begin{equation*}
B\left( \alpha \right) <0\iff \alpha <\beta \iff 0<A\left( \beta \right)
\end{equation*}%
because $B\left( \beta \right) =0$ and $A\left( \alpha \right) =0$.

Furthermore $B^{\prime }\left( \beta \right) >0$ (since $B$ is increasing in
this regime) and $\beta >0$. Thus%
\begin{equation*}
A\left( \beta \right) >0\iff \frac{A\left( \beta \right) }{\beta B^{\prime
}\left( \beta \right) }>0
\end{equation*}
Using this equivalence finally leads
to the main result in~\eqref{eq:rs}.

\section{Qualitative analysis}

\label{sct:analysis}

\noindent This section applies the main result (Theorem~\ref%
{thm:rs}) in several settings to gain some
insights into the performance of PCA.
\subsection{Dependence on balance of noise variances}

\noindent Here we would like to understand how the balance of the noise variances
affects the performance. We consider a case with two noise levels
where we sweep the noise variances $\sigma_1^2$ and $\sigma_2^2$ while
holding fixed the average noise variance 
\begin{equation*}
\bar{\sigma}^2=p_1 \sigma_1^2 + p_2\sigma_2^2.
\end{equation*}
In particular, consider 
\begin{align*}
c&=10 & p&=\left(0.7,0.3\right) \\
\theta&=1 & \bar{\sigma}&=1.3
\end{align*}
where we sweep over $\lambda \in [0,1]$ and set 
\begin{align*}
\sigma_1^2 &= \frac{\lambda}{p_1\lambda+p_2(1-\lambda)} \bar{\sigma}^2 \\
\sigma_2^2 &= \frac{(1-\lambda)}{p_1\lambda+p_2(1-\lambda)} \bar{\sigma}^2 .
\end{align*}
As intended, this fixes the average noise variance at 
\begin{equation*}
\frac{1}{n}\sum_{i=1}^n \eta_i^2 = p_1\sigma_1^2+p_2\sigma_2^2=\bar{\sigma}%
^2.
\end{equation*}
Sweeping over $\lambda$ adjusts the breakdown of the average noise variance $%
\bar{\sigma}^2$ across the two noise levels. 
%
It is not initially obvious whether better performance will occur halfway
when $\lambda=1/2$ or at the extremes when $\lambda=0$ or $\lambda=1$%
. When $\lambda=1/2$ both noise levels are the same and so all the samples
have noise variance $\bar{\sigma}^2$ and this reduces to the previously
considered homoscedastic case analyzed as an example in~\cite%
{benaych2012tsv}.
When $\lambda=0$ or $\lambda=1$, some of the points have no noise (and so PCA may do better), but the rest have noise larger than $\bar{\sigma}^2$ (and so PCA may do worse).

\begin{figure}[t]
\hrule
\centering
\begin{subfigure}[t]{0.47\linewidth}
\includegraphics[width=\linewidth]{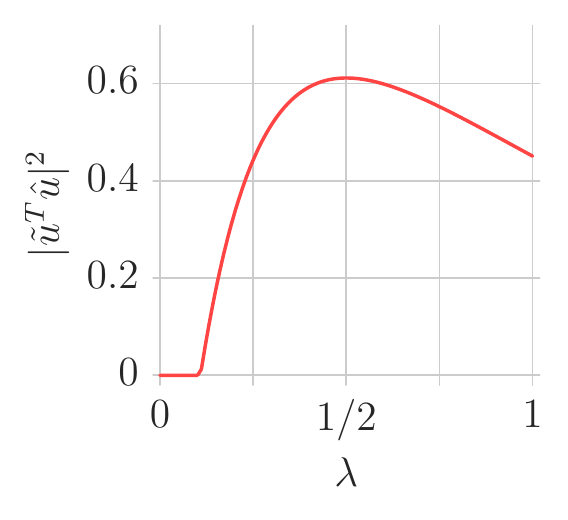}
\caption{Sweeping over noise levels while keeping the average noise variance fixed at $\bar{\sigma}^2=1.69$.
When $\lambda=1$, $\sigma_1^2$ is large and $\sigma_2^2=0$. When $\lambda=0$, $\sigma_1^2 = 0 $ and $\sigma_2^2$ is even larger.}
\label{fig:qual1a}
\end{subfigure}
\hfill
\begin{subfigure}[t]{0.47\linewidth}
\includegraphics[width=\linewidth]{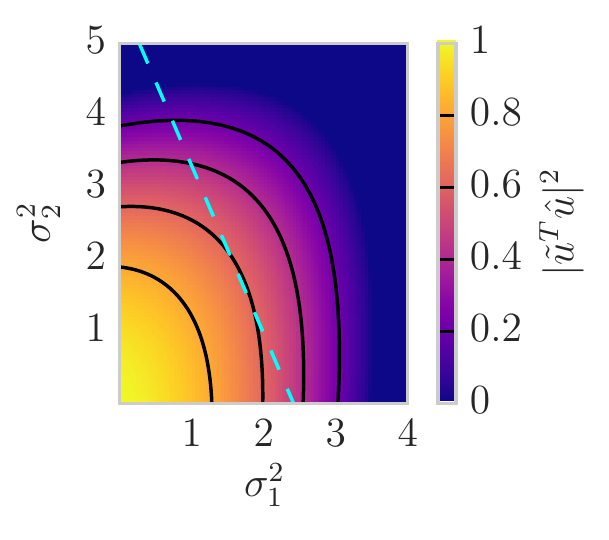}
\caption{Sweeping over both noise levels independently. The solid black curves are contours. On the dotted cyan line, the average noise variance is $\bar{\sigma}^2=1.69$.}
\label{fig:qual1}
\end{subfigure}
\caption{Asymptotic prediction under various noise levels for
$c=10$, $p=\left(0.7,0.3\right)$, $\theta=1$.
Namely, $70\%$ of samples have noise variance $\sigma_1^2$ and $30\%$ have noise variance $\sigma_2^2$.}
\hrule
\end{figure}

Figure~\ref{fig:qual1a} shows that the asymptotic prediction has a peak at 
$\lambda = 1/2$; recovery is best when the two noise levels are
the same. In other words, having samples with more noise hurts more than having samples with
corresponding less noise helps, regardless of which set is larger.
This seems to be a general phenomenon; the same has occurred for other choices of parameters we tried.

To further investigate, we use the same parameters $c,p,\theta$ as before
but sweep over both $\sigma_1^2$ and $\sigma_2^2$ (independently) and
produce a heatmap of the asymptotic prediction, shown in Figure~\ref%
{fig:qual1}.
On this figure, adjusting $\lambda$ corresponds to
moving along the line (shown as light blue dashes):
\begin{equation*}
\sigma_2^2 = \frac{1}{p_2}\bar{\sigma}^2 - \frac{p_1}{p_2}\sigma_1^2,
\end{equation*}
which has slope $-p_1/p_2 = -7/3$. Along this line, the
prediction does indeed decrease away from the diagonal $\sigma_1^2=%
\sigma_2^2$.

\begin{figure}[t!]
\hrule
\begin{center}
\vspace{2mm}
\includegraphics[width=0.81\linewidth]{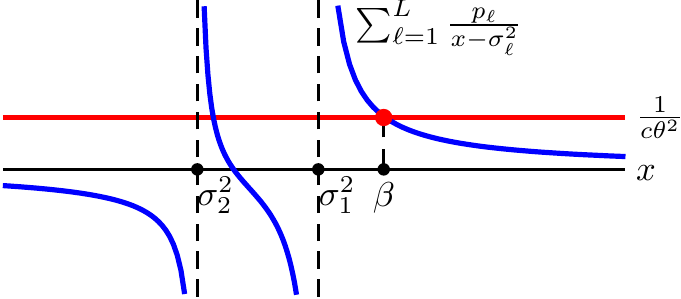}
\end{center}
\caption{Location of the largest root $\protect\beta$ of $B(x)$ for $\protect%
\lambda = 0.7$ (i.e., $\protect\sigma_1^2 = 2.04$ and $\protect\sigma_2^2 =
0.874$), $c = 10$, $p = (0.7,0.3)$ and $\protect\theta = 1$.}
\label{fig:qualB}
\hrule
\end{figure}

Figure~\ref{fig:qual1} illustrates that the prediction seems to depend primarily on the larger of the two noise
variances. This is initially surprising but might be understood by
considering the root $\beta $ of $B$. Recall that it is the largest value
$x$ satisfying%
\begin{equation*}
\frac{1}{c\theta ^{2}}=\sum_{\ell =1}^{L}\frac{p_{\ell }}{x-\sigma _{\ell
}^{2}}
\end{equation*}%
as illustrated in Figure~\ref{fig:qualB}. This figure suggests that the
largest root is heavily influenced by the largest noise variance (it is the
nearest pole). As a result, changing the other noise variances has much less impact on $\beta $ and on the
prediction.
The precise relative impact does depend on the proportions, as seen in Figure~%
\ref{fig:qual1}, where the shape of the level curves are not symmetric around 
$\sigma_1^2=\sigma_2^2$ (if the performance depended exclusively on the max
of $\sigma_1^2$ and $\sigma_2^2$, it would necessarily be symmetric).
Nevertheless, for any proportion, large noise variances can drown out the influence of small noise variances.

This insight gives a rough explanation of why it may be generally preferable to have equal noise variances ($\lambda = 1/2$) for some fixed average noise variance. Imbalance ($%
\lambda \neq 1/2$) means that one of the noise variances will be larger and
cause the performance to decline even though the other noise variance
is smaller.

\subsection{Dependence on sample-to-dimension ratio and average noise
variance}

\noindent Now consider holding everything constant except for the sample-to-dimension
ratio and average noise variance. We first suppose that there is only one
noise level (alternatively, two noise levels that are equal). In particular
we consider 
\begin{align*}
\theta&=1 & p_1 &= 1 & \ \sigma_1^2&=\bar{\sigma}^2
\end{align*}
and we sweep over $c > 1$ and $\bar{\sigma}^2 > 0$, as shown in Figure~\ref%
{fig:qual2a}. Note that this is the homoscedastic case analyzed as an
example in~\cite{benaych2012tsv} and Figure~\ref%
{fig:qual2a} illustrates the predicted phase
transition at $c=\bar{\sigma}^4$. In fact,
the asymptotic prediction in~\eqref{eq:rs} specializes to the
prediction in~\cite{benaych2012tsv} for the case where there is only one
noise level and hence the noise is homoscedastic.

\begin{figure}[t]
\hrule
\centering
\begin{subfigure}[t]{0.47\linewidth}
\includegraphics[width=\linewidth]{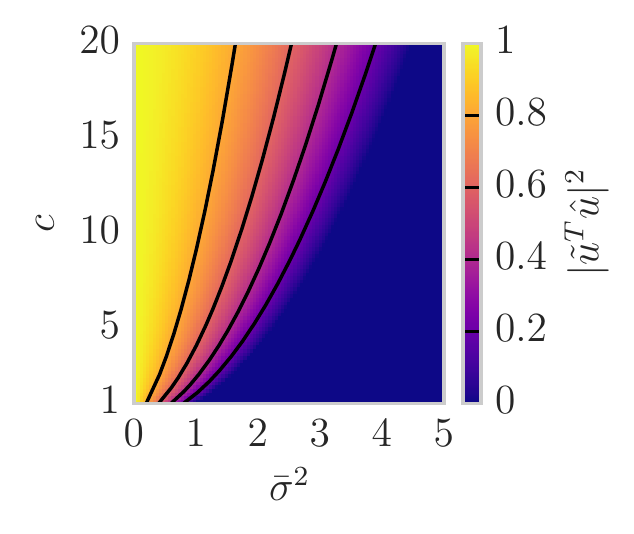}
\caption{Homoscedastic (i.e., identically distributed) noise.}
\label{fig:qual2a}
\end{subfigure}
\hfill
\begin{subfigure}[t]{0.47\linewidth}
\includegraphics[width=\linewidth]{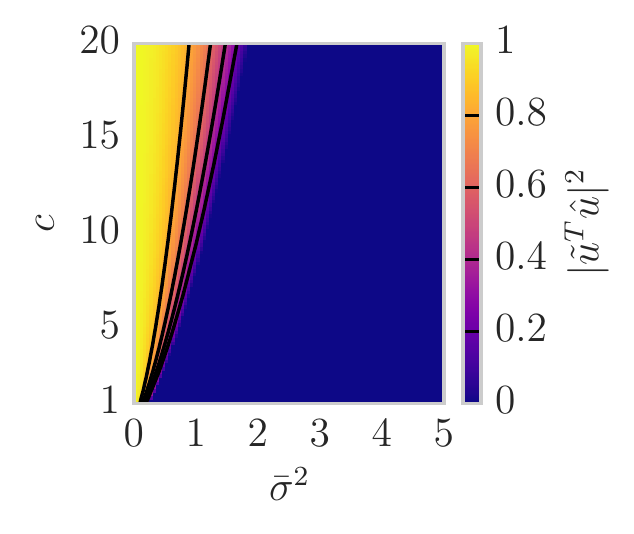}
\caption{Heteroscedastic (i.e., imbalanced) noise.}
\label{fig:qual2b}
\end{subfigure}
\caption{Asymptotic prediction as a function of average noise variance $\bar{\protect\sigma}^2$ and sample-to-dimension ratio $c$. Contours are overlaid in black. Note that the phase transition in (b) is further left than in (a); more samples are needed to tolerate the same amount of noise.}
\hrule
\end{figure}

We now consider an analogous setting where the noise is imbalanced
(i.e., heteroscedastic) with 
\begin{align*}
p_1 &= 0.9 & \sigma_1^2 = \frac{1}{2p_1} \bar{\sigma}^2 \\
p_2 &= 0.1 & \sigma_2^2 = \frac{1}{2p_2} \bar{\sigma}^2
\end{align*}
so that 
\begin{equation*}
p_1\sigma_1^2+p_2\sigma_2^2 = \frac{1}{2}\bar{\sigma}^2 + \frac{1}{2}\bar{%
\sigma}^2 = \bar{\sigma}^2.
\end{equation*}
Figure~\ref{fig:qual2b} illustrates a similar behavior with a phase transition further left. Namely, more samples are needed than in the
homoscedastic setting for the same average noise variance. This agrees
with the previous observation that performance for a given average noise
variance is best when all the points have the same noise variance (i.e., are homoscedastic).

\addtolength{\textheight}{-4.4cm}

\subsection{Dependence on sample proportions}

\noindent Finally we revisit the sweep carried out in the numerical experiments of Section~\ref{sct:experiment}. Recall that everything but the proportions were fixed. In particular,
\begin{align*}
c &= 10 & \theta&=1 & \sigma&=(1.8,0.2)
\end{align*}
and $p_2$ varied from $0$ to $1$ with $p_1 = 1-p_2$.
Figure~\ref{fig:exp} shows the prediction as a red curve (identical in both sub-figures).

As expected, the performance is best when $p_2 =
1 $ and all the samples have the lower noise variance; it is preferable to have a larger proportion of low noise samples. Interestingly, the
benefit of having more low noise samples
is not uniform through the range. The slope for small $p_2$ (close to zero)
is less steep than that for high $p_2$ (close to one). Hence, having a larger proportion of low noise samples is not as helpful when there are only a few low noise samples otherwise. A careful investigation of this phenomenon would be an
interesting area of further study.


\section{Discussion and extensions}

\label{sct:discussion}

\noindent This paper considered PCA when noise is heteroscedastic and provided a step towards the analysis of the
recovery of the subspace. In particular, we provided a simple asymptotic
prediction for the recovery of a one-dimensional subspace by PCA from noisy
heteroscedastic samples. We provided an example, illustrating
how the simple form enables easy and efficient calculation of the asymptotic
prediction, as well as an experimental verification of the prediction in
simulation. Next, we used the simple form to reason qualitatively about the
asymptotic prediction and gain new insights about the performance of PCA.
Namely, we found that the performance seems to often be most heavily
influenced by the largest noise variance present in the data. Hence,
heteroscedasticity tends to have a negative impact on the performance of PCA.

There are many avenues for potential extensions and further work. A natural
direction is to extend this work to multi-dimensional subspaces. Another
avenue of future work will be to consider a weighted version of PCA, where
the samples are first weighted in the objective function~\eqref{eq:cost} to reduce the impact of very noisy points.
Unfortunately, applying these weights violates the construction of the
subspace coefficients as being identically distributed and so this is a challenging extension. Other avenues include further investigation of
the phenomena discussed in the qualitative studies above as well as further
study of the algebraic structure of the expressions in the prediction.

\bibliographystyle{IEEEtran}
\bibliography{IEEEabrv,refs}

\end{document}